\documentclass{article}

\usepackage[utf8]{inputenc}
\usepackage{amsmath, amsfonts, dsfont}
\usepackage{graphicx}

\usepackage[numbers,square]{natbib}

\usepackage{geometry}
\title{Modeling multi-stage decision optimization problems}

\author{
Ronald Hochreiter
}

\date{April 2014}

\begin{document}

\maketitle

\begin{abstract}
Multi-stage optimization under uncertainty techniques can be used to
solve long-term management problems. Although many optimization modeling
language extensions as well as computational environments have been
proposed, the acceptance of this technique is generally low, due to the
inherent complexity of the modeling and solution process. In this paper
a simplification to annotate multi-stage decision problems under
uncertainty is presented - this simplification contrasts with the common
approach to create an extension on top of an existing optimization
modeling language. This leads to the definition of meta models, which
can be instanced in various programming languages. An example using the
statistical computing language R is shown.
\end{abstract}

\section{Introduction}\label{introduction}

We consider a multi-stage stochastic decision optimization framework
based on a discrete-time decision process, i.e.~there is a sequence of
decisions at decision stages $t = 0, \ldots, T$ where at each stage $t$
a decision taker observes the realization of a random variable $\xi_t$,
and takes a decision $x_t$ based on all observed values
$\xi_0, \ldots, \xi_t$. At the terminal stage $T$ a sequence of
decisions $x = (x_0, \ldots, x_T)$ with respective realizations
$\xi = (\xi_0, \ldots, \xi_T)$ leads to some cost $f(x, \xi)$. The goal
is to find a sequence of decisions $x(\xi)$, which minimizes a
functional of the cost $f(x(\xi), \xi)$. Multi-stage means that there is
at least one intermediary stage between root stage and terminal stage.

The design goal of the approach presented in this paper is to design a
modeling language independent of (a) the optimization modeling approach,
e.g.~expectation-based convex multi-stage stochastic programming or
worst-case optimization, as well as (b) the underlying solution
technique, e.g.~either solving a scenario tree-based deterministic
equivalent formulation or computing upper and lower bounds using primal
and dual decision rules. Finally the modeling language should (c) be
completely independent from a concrete programming language (C/C++, R,
MatLab, Python, \ldots). The idea is to compose a meta model and
instance concrete implementations semi-automatically.

\begin{figure}
\begin{center}
\scalebox{0.5}{\includegraphics{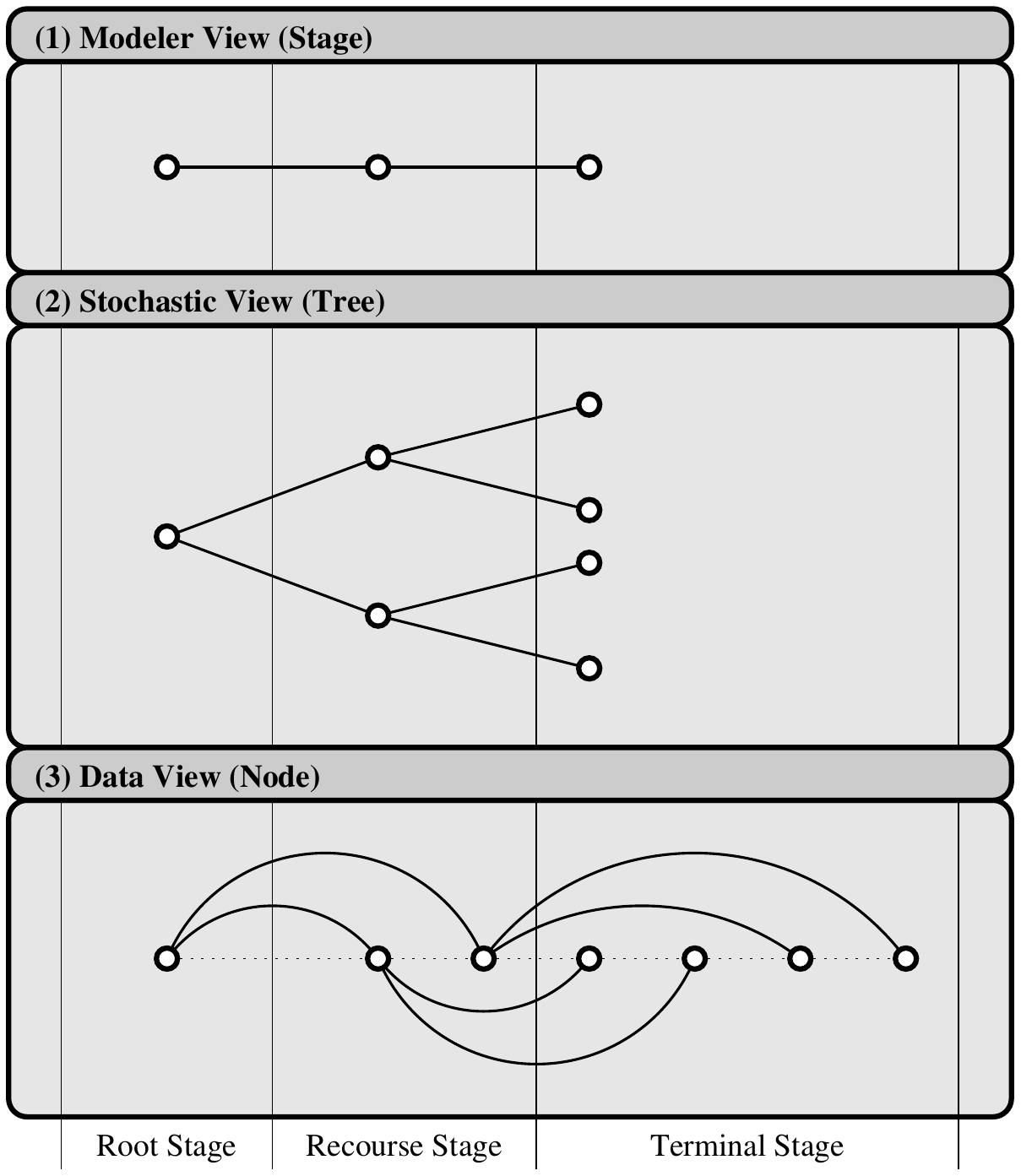}}
\end{center}
\caption{Scenario tree-based three-layered approach.}
\label{fig:msm1}
\end{figure}

Consider the two most common ways to solve multi-stage decision
optimization problems, which is on one hand the scenario-based
three-layered approach as shown in Fig. \ref{fig:msm1}. See
\cite{RuszczynskiShapiro2003} for an overview of the area of stochastic
programming, and \cite{WallaceZiemba2005} for stochastic programming
languages, environments, and applications. More information on the
modeling aspect can be found in \cite{KingWallace2013}.

\begin{figure}
\begin{center}
\scalebox{0.5}{\includegraphics{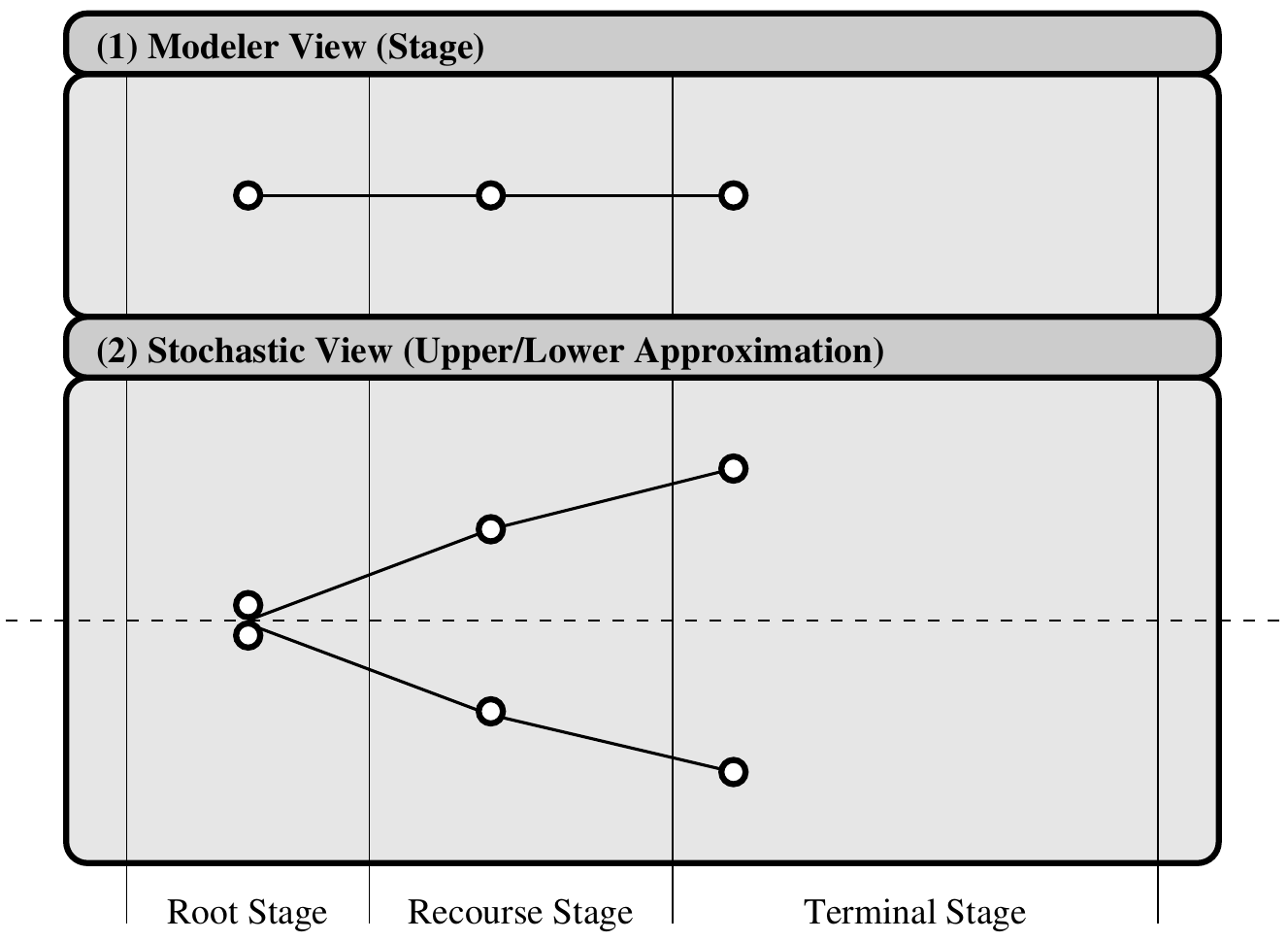}}
\end{center}
\caption{Scenario tree-free approximation.}
\label{fig:msm2}
\end{figure}

The same decision problem may also be solved using scenario tree-free
approximations \cite{KuhnEtAl2011} as shown in Fig. \ref{fig:msm2}. The
modeling language should be flexible enough to allow for applying any
solution method, i.e.~not being based on scenario trees, which is what
most modeling language extensions for multi-stage models are proposing,
see e.g. \cite{ThenieEtAl2007}, \cite{ValenteEtAl07}, \cite{Lopes2003},
\cite{DelftVial2004}, and \cite{colombo2009structure}.

\section{Multi-stage Convex Stochastic
Programming}\label{multi-stage-convex-stochastic-programming}

Consider a multi-variate, multi-stage stochastic process $\xi$ and a
constraint-set $\mathcal{X}$ defining a set of feasible combinations
$(x, \xi)$. The set $\mathcal{N}$ of functions $\xi \mapsto x$ are such
that $x_t$ is based on realizations up to stage $t$, i.e.~only
$(\xi_0, \ldots, \xi_t)$. These are the non-anticipativity constraints.
This leads to the general formulation shown in Equ. (\ref{eq:ms1}).

\begin{eqnarray}
\begin{array}{lll}
\mbox{minimize } x: & \mathds{F} \big ( f(x(\xi), \xi) \big ) & \\
\mbox{subject to} & (x(\xi), \xi) \in \mathcal{X} & \\
& x \in \mathcal{N} &
\end{array}
\label{eq:ms1}
\end{eqnarray}

The most common way to solve such a problem is to create a scenario tree
approximation of the underlying stochastic process and to build a
deterministic equivalent formulation. The problem is that most modeling
environments and languages are solely focussing on this type and mostly
provide linear-only models due to solvability concerns. Furthermore,
most allow for text-book applications only. There is almost no
flexibility provided to extend models to use real-world objective
functions and constraints.

The proposed solution is based on a complete decoupling of any scenario
tree type of modeling from the decision problem modeling process, as
shown in Fig. \ref{fig:msm1}. On the decision problem (modeling) layer
one should only be concerned with actions and decisions at stages. Other
layers differ depending on the chosen solution method. In case of
scenario trees and deterministic equivalent formulations there is an
explicit decoupling of modeling and (scenario) tree handling, i.e.~a
scenario tree layer, whose focus is to create a scenario tree which
optimally represents the subjective beliefs of the decision taker at
each node. Furthermore there is an additional data layer, which handles
the way how to (memory-)optimally store large scenario trees, access
ancestor tree nodes quickly, and other computational (tree) operations.

\section{Multi-stage Modeling
Example}\label{multi-stage-modeling-example}

Consider the stylized simple multi-stage stochastic programming example
from \cite{HeitschEtAl2006}, which is shown in Equ. (\ref{eq:ex1}).

\begin{eqnarray}
\begin{array}{lll}
\mbox{minimize} & \mathds{E} \big ( \sum_{t=1}^T V_t x_t \big ) & \\
\mbox{subject to} & s_t - s_{t-1} = x_t & \forall t = 2, \ldots, T \\
& s_1 = 0, s_T = a, & \\
& x_t \geq 0, s_t \geq 0. & 
\end{array}
\label{eq:ex1} 
\end{eqnarray}

The decision to be computed with this model is the optimal purchase over
time under cost uncertainty, where the uncertain prices are given by
$V_t$, and the decisions $x_t$ are amounts to be purchased at each time
period $t$. The objective function aims at minimizing expected costs
such that a prescribed amount $a$ is achieved at $T$; $s_t$ is a state
variable containing the amount held at time $t$.

In Tab. \ref{tab:ex1} a concise meta formulation of this problem can be
seen. The general syntax is borrowed from algebraic modeling languages
like AMPL \cite{FourerEtAl2002} and ZIMPL \cite{Koch2004}.

\begin{table}
\caption{Modeling formulation of Equ. (\ref{eq:ex1}).}
\begin{verbatim}
deterministic a: T;
stochastic x, s, objective_function: 0..T;
stochastic non_anticitpativity: 1..T;
stochastic root_stage: 0;
stochastic terminal_stage: T;

param a;
var x >= 0, s >= 0; 

minimize objective_function: E(V * x);
subject to non_anticitpativity: s - s(-1) = x;
subject to root_stage: s = 0;
subject to terminal_stage: s = a;
\end{verbatim}
\label{tab:ex1}
\end{table}

The most striking feature is that any relation to stages is removed from
the definition of the optimization model - parameters, variables,
objective function, and constraints. To accommodate for the definition
of stages, the proposed stochastic modeling language contains two
additional keywords for any of these objects, i.e.

\begin{itemize}
\itemsep1pt\parskip0pt\parsep0pt
\item
  \texttt{deterministic} \texttt{\it objects}: \textit{stage-set};
\item
  \texttt{stochastic} \texttt{\it objects}: \textit{stage-set};
\end{itemize}

Speaking in scenario tree notation the stochastic objects are defined on
the underlying node structure and deterministic objects are defined on
the stage structure, i.e.~the latter contain the same value for all
nodes in the respective stage. To define stochastic objective functions
and stage recourse the following functions are defined, e.g.~the most
commonly used expectation functional for objective functions is simply
expressed by the function \texttt{E()}. Furthermore, there is a special
way to define stage-wise recourse for stochastic variables, i.e.
\texttt{variable-name(recourse-depth)}. Note that while most modeling
languages allow for a single stage recourse only, this definition allows
for any number of recourse stages.

\begin{table}
\caption{Implementation of model \ref{tab:ex1} using the language R.}
\begin{verbatim}
m <- model()
parameter(m, a)
variable(m, x, lb=0)
variable(m, s, lb=0) 

minimize(m, "objective", "E(V * x)")
subject_to(m, "non_anticitpativity", "s - s(-1) = x")
subject_to(m, "root_stage", "s = 0")
subject_to(m, "terminal_stage", "s = a")

deterministic(m, "T", a)
stochastic(m, "0..T", x, s, "objective")
stochastic(m, "1..T", "non_anticipativity")
stochastic(m, "0", "root_stage")
stochastic(m, "T", "terminal_stage")

optimize(m)
\end{verbatim}
\label{tab:ex1r}
\end{table}

Tab. \ref{tab:ex1r} shows the modeling example in some concrete
implementation for the statistical computing language R \cite{R2013}.
This definition can be easily converted to a deterministic equivalent
formulation or any other reformulation - all information is available in
a concise format.

\section{Conclusion}\label{conclusion}

In this paper, a modeling language framework for a successful simplified
meta modeling of multi-stage decision problems under uncertainty is
shown, which allows for automatic reformulation and solution of
multi-stage problems. This can be seen as a basis to build a model-based
multi-stage problem library, especially because of its inherent
decoupling from the underlying optimization technique as well as the
fact that it is not bound to a specific programming language.
Furthermore it is easy to integrate robust and stochastic optimization
techniques to allow for comparing solutions to determine, which approach
is optimally suited for which class of decision models. There are many
ways to extend the proposed meta language - possible straight-forward
extensions are e.g.~quantiles for objective functions. In addition,
application-related risk measures (shortcuts) can be defined, e.g.
\texttt{CVaR(objects)}, as well as probabilistic constraints.

\bibliographystyle{plainnat}
\bibliography{multistage}

\end{document}